\documentclass[12pt,oneside]{article}
\usepackage{amsmath,amssymb}
\numberwithin{equation}{subsection}
\def\overset#1#2{{\mathop{\kern0pt#2}\limits^{#1}}}
\begin{document}
%

\begin{center}
\large\bf
THE ACTION OF CONTACT TRANSFORMATIONS PSEUDOGROUP
ON THE SECOND ORDER ODES WHICH ARE CUBIC IN SECOND DERIVATIVE
\\[5pt]
\rm Vadim V. Shurygin, jr.\footnote
{The author is supported by RFBR  grant  15-31-50220.}\\[5pt]
{\it Kazan Federal University, Russia}
\end{center}

\begin{abstract}
In the present paper we consider the problem of local
equivalence of second order ODEs which are cubic in second
derivative under the action of the pseudogroup of contact
transformations.
We show how it may be reduced to the equivalence problem
of 2-webs in $\mathbb{R}^3$ under the action of finite-dimensional group, admitted by the associated second order equation.
\end{abstract}



\subsection{Introduction}

The problem of equivalence of second order ODEs of the form
\begin{equation}
\label{du1}
y''=F(x,y,y')
\end{equation}
under the action of some transformations pseudogroup was the
subject of many papers.
S.\,Lie proved that any two such equations are equivalent under the
action of the pseudo\-group of contact transformations on the space
$\mathbb{R}^3(x,y,p)$.
A.\,Tresse~\cite{Tr2} found
the complete set of differential invariants of the ODE~(\ref{du1})
under the action of the pseudogroup of point transformations.
B.\,Kruglikov in his review paper~\cite{Krug} presented
the complete description of  the algebra of differential invariants
and the solution of  the problem of point equivalence of such ODEs.

We consider the problem of equivalence of the ODEs of the form
\begin{equation}
\label{osn-ur}
y''{}^n + A_1(x,y,y')y''{}^{n-1} + \ldots +
A_{n-1}(x,y,y')y'' + A_n(x,y,y')=0
\end{equation}
under the action of the pseudogroup of contact transformations.
We suppose that this equation has
 $n$ distinct roots
$$
y''=\lambda_i(x,y,y'),\quad i=1,\dots, n.
$$
The case $n=2$ was solved in~\cite{Vshjr12}.
The equivalence problem of such ODEs reduces to the problem of
point equivalence of some ODEs of the form (\ref{du1}),
associated with the original equations.
In the present paper we show how the case
$n=3$ may be reduced to the case $n=2$.
Moreover, we consider the most interesting situation for
$n=3$, namely, the case when the associated equation for
the quadratic equation $(y''-\lambda_1)(y''-\lambda_2)=0$,
is linearizable.
This case reduces to the equivalence problem of 2-parametric
families of integral curves of Cartan distribution in
$\mathbb{R}^3$ under the action of projective group $SL_3$ admitted by
the associated equation.
Such families were called  Legendrian 2-webs in~\cite{Wang}.
This problem can be reformulated as a classification problem of
ODEs of the form (\ref{du1}) under the action of $SL_3$.


\subsection{Basic definitions and constructions}

Let $\pi$ be a vector bundle and $G$ be a pseudogroup of
diffeomorphisms acting on $\pi$.
The action of $G$ naturally lifts to the action
on the space $J^k\pi$ of $k$-jets of sections of
$\pi$.

By a smooth function $f\in C^\infty(J^\infty\pi)$ we mean
the smooth function on the finite-dimensional
jet space $J^k(\pi)$ for some  $k\ge 1$.

\smallskip
{\bf Definition.}
The function $I\in C^\infty(J^k\pi)$ is called the
{\it (absolute) scalar differential
invariant  of $k$-th order}, if it is constant along the orbits
of the lifted action of $G$ on $J^k\pi$.
\smallskip

The space $\mathcal I$ of all invariants forms an algebra with
respect to  algebraic operations of linear combinations over $\mathbb{R}$,
multiplication and  the composition
$I_1, \dots , I_m \to I = F(I_1, \dots , I_m)$
for any smooth function $F\in C^\infty(\mathbb{R}^m)\to
\mathbb{R}$, $m\ge 1$.

Let  $G$ be a connected Lie group and let
$\mathfrak g$ be its Lie algebra.
For $X\in \mathfrak g$ denote by $\widehat X$ the lift
of  $X$ to $J^\infty(\pi)$ (see~\cite{KLV}). Then $I$ is an absolute differential invariant
if and only if  $L_{\widehat X}(I)=0$ for all $X\in \mathfrak g$.
For more information on differential invariants see, e.g.,~\cite{KLV, KL, Ol}.

\smallskip
{\bf Definition.}
The function $F\in C^\infty(J^k\pi)$
is called the {\it relative scalar differential
invariant of $k$-th order}, if for all $g\in G$ there holds
$$
g^* F=\mu(g)\cdot F ,
$$
where $\mu:G\to C^\infty(J^\infty\pi)$ is
a smooth function satisfying the condition
$$
\mu(g\cdot h)=
h^*\mu(g)\cdot \mu(h),\quad \mu(e)=1.
$$
In other words, the equation $F=0$ is invariant under the action
of $G$.
The function $\mu$ is called the {\it weight function}.
\smallskip

The infinitesimal analogue of this definition is
that $F$ is a relative invariant if and only if
$L_{\widehat X}(F)=\mu_X\cdot F$, where the map
$\mu:\mathfrak{g}\to C^\infty(J^\infty(\pi))$ satisfies the
condition
$$
\mu_{[X,Y]}=L_{\widehat X}(\mu_Y)-L_{\widehat Y}(\mu_X), \qquad
\forall\,\, X,Y\in\mathfrak{g}.
$$

Let $(x_1,\dots, x_n) $ be the coordinates in the base manifold of
the bundle $\pi$.

\smallskip
{\bf Definition.}
A $\mathfrak g$-{\it invariant derivation} is a linear combination
of total derivatives
$$
\nabla=\sum_{i=1}^n A_i \dfrac{d}{dx^i}, \quad A_i\in
C^\infty(J^\infty\pi),~i=1,\dots,n,
$$
which commutes with all prolongations of vector fields
$X\in\mathfrak{g}$.

\smallskip
For each differential invariant $I$ the function
$\nabla(I)$ also is a differential invariant (usually, of order, higher
than the order of $I$).
This fact allows one to obtain new invariants out of known ones
using invariant derivations.
For the algebra $\mathcal I$ of differential invariants
the {\it Lie-Tresse theorem} is valid.
It claims that there exists a finite number of basic differential
invariants and basic invariant derivations, such that any other
differential invariant can be expressed in terms of these basic invariants
and its derivations (see~\cite{KL, Kum, Ol}).

The following result of V.\,Lychagin~\cite{Lych08} allows one to find
invariant derivations.
Let $G$ be a Lie group acting on
$\mathbb{R}^n$.
Let us identify the elements of its
Lie algebra $\mathfrak g$  with contact vector fields
$X_f$ on $J^1(\mathbb{R}^n)$ (here $f$ is the generating function of $X_f$).

\smallskip
{\bf Proposition 1.}
{\it
Let $x_1$, \dots, $x_n$ be coordinates in $\mathbb{R}^n$,
and let $(x_1, \dots, x_n, u, p_1, \dots, p_n)$ denote
the corresponding canonical coordinates in $J^1 (\mathbb{R}^n)$.
Then the derivation
$$
\nabla=\sum\limits^{n}_{i=1}
A_i
\dfrac{d}{dx^i}
$$
is $\mathfrak g$-invariant if and only if
functions $A_i \in C^\infty(J^\infty(\mathbb{R}^n))$,
$i =1, \dots, n$, satisfy the follo\-wing system of PDEs
$$
X_f (A_i) +\sum\limits^n_{j=1}
\dfrac{d}{dx^j}\left(\dfrac{\partial f}{\partial p_i}\right) A_j = 0
$$
for all $j=1,\dots,n$, $X_f\in\mathfrak g$. }
\smallskip

Let $\mathcal{ C}$ denote the standard contact distribution
(Cartan distribution) on
$\mathbb{R}^3(x,y,p)$
and let
$\omega=dy- p\, dx$ be the contact form.
Geometrically each of the equations
$y''=\lambda_i(x,y,y')$ determines the 1-dimensional
distribution ${\mathcal{F}}_i$ in the 3-dimensional contact
manifold
$\mathbb{R}^3(x,y,p)$.
This distribution
is given by a vector field
\begin{equation}
\label{Xi}
X_i=\dfrac{\partial}{\partial x} + p
\dfrac{\partial}{\partial y } +
\lambda_i(x,y,p) \dfrac{\partial}{\partial p}.
\end{equation}
and lies in $\mathcal{C}$.

In the paper~\cite{Vshjr12} we solved the problem of contact equivalence
for the equation~(\ref{osn-ur}) when $n=2$:
\begin{equation}
\label{du2}
\mathcal{E}=\{(y''-\lambda_1(x,y,y'))(y''-\lambda_2(x,y,y'))=0\}.
\end{equation}
Let $X_1$ and $X_2$ be the vector fields~(\ref{Xi}) corresponding
to the equation  (\ref{du2}) and let $a_1$, $b_1$ be any two
independent integrals
of  $X_1$, and $a_2$, $b_2$ be independent integrals of $X_2$.
Any three of these four functions are independent and the fourth one may
be expressed in terms of these three. Let, for example,
$b_2=h(a_1,b_1,a_2)$.
Let us introduce the new coordinate system
$$
a=a_1,~~ b=b_1,~~ c=-\dfrac{h_{a_1}}{h_{b_1}}
$$
(the lower indices
$a_i$, $b_i$ denote the partial derivatives with respect to $a_i$,
$b_i$).
In these coordinates the fields $X_1$, $X_2$ have the following
form (up to a multiple)
$$
X_1=\dfrac{\partial}{\partial c},\quad
X_2=\dfrac{\partial}{\partial a} +
c\dfrac{\partial}{\partial b} + G_2(a,b,c)
\dfrac{\partial}{\partial c},
$$
where
$$
G_2(a,b,c) =- \dfrac{h_{a_1a_1}h_{b_1}^2 -2h_{a_1}h_{b_1}h_{a_1b_1} +
h_{b_1b_1}h_{a_1}^2 }{h_{b_1}^3}.
$$
The exterior 1-form $\theta=db- c\, da$ determines
the contact structure on $\mathbb{R}^3$.

We will say that the ODE
$$
\mathcal{E}^a_2=\{y''=G_2(x, y, y')\},
$$
determined by  $X_2$, is {\it associated} with $\mathcal{E}$.
The choice of another pairs of integrals of $X_1$ and $X_2$ leads
to another associated equation which is point equivalent to the
previous one.
Moreover, if we interchange  $X_1$ and $X_2$, then we obtain
another equation  $\mathcal{E}^a_1=\{y''=G_1(x, y, y')\}$ associated with
$\mathcal{E}$.
We will say that these two equations are {\it dual}.
Two dual equations are not necessary point equivalent  (see~\cite{Vshjr12}).
Thus, it is correct to speak about the two equivalence classes
$([\mathcal{E}^a_1], [\mathcal{E}^a_2])$ associated with~$\mathcal{E}$.

\smallskip
{\bf Theorem 1.~\cite{Vshjr12}}
{\it
Let
$$
\mathcal{E}=\{(y''-\lambda_1)(y''-\lambda_2)=0\}
\hbox{~~~and~~~}
\widetilde{\mathcal{E}}=\{(y''-\widetilde\lambda_1)(y''-\widetilde\lambda_2)=0\}
$$
be two ODEs of the form {\rm (\ref{du2})}
and let $([\mathcal{E}^a_1], [\mathcal{E}^a_2])$ and
$([\widetilde{\mathcal{E}}^a_1], [\widetilde{\mathcal{E}}^a_2])$
be the pairs of  equivalence classes of their associated equations.
Then $\mathcal{E}$ and $\widetilde{\mathcal{E}}$
are equivalent under the action of the pseudo\-group of contact trans\-for\-mations
if and only if
one of the classes
$([\mathcal{E}^a_1], [\mathcal{E}^a_2])$
coincides with one of the classes
$([\widetilde{\mathcal{E}}^a_1], [\widetilde{\mathcal{E}}^a_2])$.}

\subsection{From ODEs to foliations}

It follows that every equation $y''=\lambda(x,y,y')$
is equivalent to the
foliation ${\cal F}_\lambda$ formed by the curves with equations
\begin{equation}
\label{abc}
a(x,y,p)={\rm const}, \quad b(x,y,p)={\rm const}
\end{equation}
in $\mathbb{R}^3(x,y,p)$.
The functions $a$ and $b$ satisfy the conditions
$$
da\wedge db \wedge \omega=0, \quad da\wedge db \ne0,
$$
where $\omega$ is the Cartan form.
These functions are essential up to
diffeomorhisms
\begin{equation}
\label{ab-z}
(a,b)\to
(\varphi(a,b), \psi(a,b)).
\end{equation}

Recall that a {\it $k$-web}
in a manifold $M$ is a set of $k$ one-dimensional foliations in $M$ which
are in general position, i.e., the tangent vector fields
corresponding to different foliations
are linear independent.
For a contact structure $\omega$ on $M$, a Legendrian curve is a curve $\gamma$
such that $\omega|_\gamma=0$, i.e., the curve which is tangent to the contact distribution.
A {\it Legendrian $k$-web} is  a set of $k$ pairwise transversal foliations by
Legendrian curves (see J.S.\,Wang~\cite{Wang}).
The foliation ${\cal F}_\lambda$ is a Legendrian 1-web.
It follows that
every ODE (\ref{osn-ur}) is equivalent to a Legendrian $n$-web.

Let now the equation
\begin{equation}
\label{du3}
y''{}^3+A(x,y,y') y''{}^2+B(x,y,y')y''+C(x,y,y')=0
\end{equation}
be given and let
$$
y''=\lambda_i(x,y,y'),\quad i=1,2,3,
$$
be its roots.
Choose any two of them, say $\lambda_1$ and $\lambda_2$,
and
let
\begin{equation}
\label{yG}
y''=G(x,y,y')
\end{equation}
be the equation associated with $(y''-\lambda_1)(y''-\lambda_2)=0$.
In the paper~\cite{Vshjr12} we showed that
the action of contact transformations pseudogroup on the
ODE (\ref{du2}) induces the action of point transformations
pseudogroup on the associated ODEs.
It follows that the group $\mathcal{G}$ of point transformations, admitted
by the equation (\ref{yG}) acts on the set of integrals of the third equation
$y''=\lambda_3(x,y,y')$.
It is well-known~(see, e.g., \cite{Lie}) that such a group may have dimension 0, 1, 2, 3
or 8.
Thus, the problem of contact equivalence
of ODEs  (\ref{du3}) reduces to the problem
of equivalence of foliations
(\ref{abc})
in $\mathbb{R}^3(x,y,p)$
under the action of the group $\mathcal{G}$.

Let  ${\cal F}_\lambda$ be a foliation (\ref{abc}).
Such a  family is a section of a 2-dimensional bundle
\begin{equation}
\label{R53}
\mathbb{R}^5\to \mathbb{R}^3,\qquad
(x,y,p,a,b)\to (x,y,p).
\end{equation}
Since the functions $a$ and $b$
are significant up to fiberwise diffeomeorphisms (\ref{ab-z}),
we pass to another bundle
$$
\pi:\mathbb{R}^5\to \mathbb{R}^3,
\qquad
(x,y,p,f,g)\to (x,y,p).
$$
The fibers of $\pi$ consist of the values of the functions
\begin{equation}
\label{fg}
f:=\dfrac{a_xb_y-a_yb_x}{a_yb_p-a_pb_y},\qquad
g:=\dfrac{a_pb_x-a_xb_p}{a_yb_p-a_pb_y}.
\end{equation}
The group of fiberwise diffeomorphisms of the bundle (\ref{R53})
acts of $\pi$ trivially.
The functions $a$ and $b$ are  not arbitrary,
they are  integrals of the  vector field
$$
X=\dfrac{\partial}{\partial x} + p
\dfrac{\partial}{\partial y } +
G(x,y,p) \dfrac{\partial}{\partial p}.
$$
It follows by direct computations that the functions (\ref{fg})
satisfy the relations
\begin{equation}
\label{fGgp}
f(x,y,p)=G(x,y,p),\quad g(x,y,p)=p.
\end{equation}
The conditions (\ref{fGgp}) determine
the 1-dimensional subbundle $\widetilde \pi$ of $\pi$.

To solve the equivalence problem of foliations (\ref{abc}) under the
action of $\cal G$ we need to find the algebra of differential invariants of
this action on $\widetilde{\pi}$.
The computations of differential invariants were performed
using the {\tt Maple} packages {\tt DifferentialGeometry} and {\tt JetCalculus}
by I.M.\,Anderson.

To do this, we need to
lift the action of  the Lie algebra $\mathfrak g$ of the Lie group $\cal G$
to the bundle $\widetilde \pi$.
We prolong
the vector fields from $\mathfrak g$ to
$\mathbb{R}^3(x,y,p)$ using the standard formulas~\cite{Ibr92, KLV, Ol2}
$$
\xi\dfrac{\partial}{\partial x}+
\eta\dfrac{\partial}{\partial y} \longmapsto
\xi\dfrac{\partial}{\partial x}+
\eta\dfrac{\partial}{\partial y}+
\zeta\dfrac{\partial}{\partial p}, \qquad
\zeta=\eta_x+p(\eta_y-\xi_x)-p^2\xi_y.
$$
The lift of the field
\begin{equation}
\label{Xpi}
X=\xi(x,y,p)\dfrac{\partial}{\partial x}+
\eta(x,y,p)\dfrac{\partial}{\partial y}+
\zeta(x,y,p)\dfrac{\partial}{\partial p}
\end{equation}
to $\pi$
has the form
\begin{multline}
\widehat{X}^{(0)}=\xi\dfrac{\partial}{\partial x}+\eta\dfrac{\partial}{\partial y}+
\zeta\dfrac{\partial}{\partial p}
+\Bigl(\zeta_p f + \zeta_y g + \zeta_x- (\xi_p f + \xi_y g +\xi_x)f\Bigr)
\dfrac{\partial}{\partial f} \,+ \\
+ \Bigl(\eta_p f + \eta_y g + \eta_x - (\xi_p f + \xi_y g+\xi_x)g\Bigr)
\dfrac{\partial}{\partial g}
\end{multline}
and its restriction
to $\widetilde \pi$ is
$$
\widehat{X}^{(0)}_{\widetilde \pi}=\xi\dfrac{\partial}{\partial x}+\eta\dfrac{\partial}{\partial y}+
\zeta\dfrac{\partial}{\partial p}
+\Bigl(\zeta_p f + \zeta_y p + \zeta_x- (\xi_p f + \xi_y p +\xi_x)f\Bigr)
\dfrac{\partial}{\partial f}.
$$

\subsection{Projective classification of
foliations}

In the present section we consider the most interesting case of the projective group
$SL(3)$ of dimension 8.
The ODEs that admit this group, are linearizable,
that is, point equivalent to the ODE
$y''=0$ (see, e.g.,~\cite{Ibr09,  Lie}).
The conditions under which the ODE is linearizable can be found
in~\cite{Ibr92, Ibr09, Krug, Mah, Tr2}.

In what follows we  assume that
the coordinate system  $(x,y,p)$ in  $\mathbb{R}^3$
is chosen in such a way that the associated equation (\ref{yG})
has the form
$y''=0$.
The Lie algebra  $\mathfrak{sl}(3)$
of the group $SL(3)$ consists of the vector fields
\begin{equation}
\label{Xsl3}
(a_0 + a_{11}x + a_{12}y)\dfrac{\partial}{\partial x}
+ (b_0 + a_{21}x + a_{22}y)\dfrac{\partial}{\partial y}
+ (c_1x+c_2y) \left( x\dfrac{\partial}{\partial x} +
y\dfrac{\partial}{\partial y}\right),
\end{equation}
on $\mathbb{R}^2(x,y)$, where $a_0$, $b_0$, $a_{ij}$, $c_i$ are constants.

\smallskip
{\bf Proposition 2.}
{\it The action of $SL(3)$ on $\widetilde{\pi}$ has the following relative invariants:
one relative invariant of order~0, one of order 1
$$
I_0=f,
\qquad
I_1=pf_yf_p-3ff_y+f_xf_p
$$
and five of order 2
$$
\begin{array}{l}
H_1=3f_{pp}f^2-2ff_p^2,
\\[5pt]
H_2=3f(f_{xx}+2pf_{xy}+p^2f_{yy})-4(pf_y+f_x)^2,
\\[5pt]
H_3=3ff_pf_{xp}+ 3f(pf_p-3f)f_{yp} + 3f(pf_y+f_x)f_{pp} +
f_p(9ff_y-5f_p(pf_y+f_x)),
\\[5pt]
H_4=3f(f_pf_{xx} + (2pf_p-3f)f_{xy} + (pf_y+f_x)(f_{xp} + pf_{yp})
+p(pf_p-3f)f_{yy}) -\\[3pt]
\hskip8cm-(pf_y+f_x)(7f_p(pf_y+f_x)-12ff_y),
\\[5pt]
H_5=3f\Bigl(f_p^2f_{xx} + 2f_p(pf_p-3f)f_{xy}+(pf_p-3f)^2f_{yy}
+2(2pf_yf_p+2f_xf_p-3ff_y)f_{xp}\, + \\[3pt]
\qquad\qquad
+\,2(2p^2f_yf_p+2pf_xf_p-6pff_y-3ff_x)f_{yp}
+(pf_y+f_x)^2f_{pp} \Bigr)\,-\\[3pt]
\hskip6cm
-\, 18f_p^2(pf_y+f_x)^2 + 60ff_yf_p(pf_y+f_x) -36f^2f_y^2.
%
%
\end{array}
$$
}

%
%

\smallskip
The corresponding weights are
$$
\begin{array}{l}
\mu(I_0)=-2a_{11}+a_{22}-3a_{12}p-3c_1x-3c_2xp,\\
\mu(I_1)=-4a_{11}+a_{22}-5a_{12}p-7c_1x-c_2(5xp+2y),\\
\mu(H_1)=-4a_{11}+a_{22}-5a_{12}p-7c_1x-c_2(5xp+2y)=\mu(I_1),\\
\mu(H_2)=2(-3a_{11}+a_{22}-4a_{12}p-5c_1x-c_2(4xp+y)),\\
\mu(H_3)=-5a_{11}+a_{22}-6a_{12}p-9c_1x-3c_2(2xp+y),\\
\mu(H_4)=-7a_{11}+2a_{22}-9a_{12}p-12c_1x-3c_2(3xp+y),\\
\mu(H_5)=2(-4a_{11}+a_{22}-5a_{12}p-7c_1x-c_2(5xp+2y))=2\mu(I_1).
\end{array}
$$

\smallskip
{\bf Definition.}
By an {\it (absolute) differential invariant} of the action of
$SL(3)$ on $\widetilde{\pi}$ we understand a function
$I\in J^\infty(\widetilde{\pi})$
which is constant along the
orbits of the action of $SL(3)$
and is polynomial in derivatives of second
order and higher and in the functions  $|I_0|^{-1/2}$ and
$|I_1|^{-1/2}$.
%

\smallskip

Since $\dim J^1(\widetilde{\pi})=7$ and $\mu(I_0)$ is not proportional to
$\mu(I_1)$, there are no differential
invariants of order less than 2.

One can verify that two weights, $\mu(I_0)$ and $\mu(I_1)$,
are generators of the linear space of weights over $\mathbb{Q}$.
Expressing the other weights in terms of these two, we get the following

\smallskip
{\bf Theorem  2.}
{\it The action of $SL(3)$ on $\widetilde{\pi}$
has five basic absolute differential invariants of order 2:}
$$
M_1=\dfrac{H_1}{I_1},
\quad
M_2=\dfrac{H_2}{I_0I_1},
\quad
M_3=\dfrac{{|I_0|}^{1/2}\cdot H_3}{|I_1|^{3/2}},
\quad
M_4
=\dfrac{H_4}{|I_0|^{1/2}\cdot |I_1|^{3/2}},
\quad
M_5=\dfrac{H_5}{I_1^2}.
$$

\smallskip

Let
$\dfrac{d}{dx}$, $\dfrac{d}{dy}$ and $\dfrac{d}{dp}$
denote the total derivatives.

\smallskip
{\bf Definition.}
By an
$\mathfrak{sl}(3)$-{\it invariant derivation} we mean the linear
combination of total derivatives
$$
\nabla=A\dfrac{d}{dx} + B\dfrac{d}{dy} + C\dfrac{d}{dp}, \qquad A,B,C\in
C^\infty(J^\infty(\widetilde{\pi})),
$$
invariant under the prolonged action of $SL(3)$. The functions
$A$, $B$, $C$ are supposed to be polynomial in derivatives of second
order and higher and in the functions $|I_0|^{-1/2}$ and
$|I_1|^{-1/2}$.

\smallskip
Using the Proposition 1, we find three basic $\mathfrak{sl}(3)$-invariant
derivations:
$$
\nabla_1 = C_1\dfrac{d}{dp}, \qquad
\nabla_2 =
A_2\dfrac{d}{dx}+B_2\dfrac{d}{dy},\qquad
\nabla_3 =
A_3\dfrac{d}{dx}+B_3\dfrac{d}{dy}+C_3\dfrac{d}{dp},
$$
 where
$$
\begin{array}{l}
C_1=\dfrac{f^{3/2}}{(pf_yf_p-3ff_y+f_xf_p)^{1/2}}=
\dfrac{|I_0|^{3/2}}{|I_1|^{1/2}},\\[15pt]
A_2= \dfrac{f^{1/2}}{(pf_yf_p-3ff_y+f_xf_p)^{1/2}}
=\dfrac{|I_0|^{1/2}}{|I_1|^{1/2}},\quad
B_2=pA_2,\\[15pt]
A_3=\dfrac{ff_p}{pf_yf_p-3ff_y+f_xf_p}=\dfrac{f_pI_0}{I_1},\quad
B_3=\dfrac{f(pf_p-3f)}{pf_yf_p-3ff_y+f_xf_p}=\dfrac{(pf_p-3f)I_0}{I_1},\\[15pt]
\qquad\qquad C_3=\dfrac{f(pf_y+f_x)}{pf_yf_p-3ff_y+f_xf_p}=\dfrac{(pf_y+f_x)I_0}{I_1}.
\end{array}
$$

The following commutation relations hold
\begin{equation}
\label{fn123}
\begin{array}{l}
[\nabla_1,\nabla_2]= \dfrac16M_4\cdot\nabla_1
-\dfrac16M_3\cdot\nabla_2
-\dfrac13\cdot\nabla_3,\\[11pt]
[\nabla_1,\nabla_3]=\dfrac16(M_5-4)\cdot\nabla_1 +
\dfrac13M_1\cdot\nabla_2
-\dfrac13M_3\cdot\nabla_3,\\[11pt]
[\nabla_2,\nabla_3]= \dfrac13M_2\cdot\nabla_1 +
\dfrac16(M_5+4)\cdot\nabla_2
-\dfrac13M_4\cdot\nabla_3.
\end{array}
\end{equation}

\smallskip
{\bf Definition.}
We call the point $x_k\in J^k(\widetilde{\pi})$ {\it regular}, if
$I_0I_1\ne0$ at this point.
This is equivalent to the fact that the $SL(3)$-orbit
passing through the regular point has  dimension equal to
$8$.
In what follows we consider only
the orbits of regular points.

\smallskip
{\bf Theorem 3.}
{\it
Algebra of differential invariants of the
action of  $SL(3)$ on $\widetilde{\pi}$
is generated by five invariants $M_1$, $M_2$,
$M_3$, $M_4$, $M_5$ and three invariant derivations $\nabla_1$,
$\nabla_2$, $\nabla_3$.
This algebra separates  regular orbits.}

\smallskip
{\bf Proof.}
There are no differential
invariants of order less than 2.
The codimension of a regular orbit $\mathcal{O}_{2}\subset J^2(\widetilde{ \pi})$
equals 5.
At the same time, there exist five independent
invariants of order 2.
It can be checked by direct computations that
at a regular point $x_2\in \dim
J^2(\widetilde{\pi})$ one has $h\cdot x_2=x_2$, $h\in SL_3$ if and only if
$h=e$. Consequently, invariants $M_1$, $M_2$, $M_3$, $M_4$, $M_5$
generate the space if invariants of second order and separate regular
orbits.

Note that invariants $M_1$, \dots, $M_5$ are linear in second order
derivatives $f_{xx}$, \dots, $f_{pp}$,
and that  $C_1$, $A_2$, $B_2$, $A_3$, $B_3$, $C_3$
are functions on $J^1(\widetilde{\pi})$.
It follows that the functions $M_{ik}:=\nabla_i M_k$,
$i=1,2,3$,
$k=1,\dots,5$, are linear
in
third order derivatives.

The fibers ${\cal F}_3$ of the bundle $J^3(\widetilde{\pi})\to J^2(\widetilde{\pi})$
are 10-dimensional.
The following 10 invariants of third order have independent
symbols
\begin{equation}
\label{nM}
M_{11}, ~ M_{12}, ~ M_{13}, ~ M_{14}, ~ M_{15}, ~ M_{21},
~ M_{22}, ~ M_{24}, ~ M_{25}, ~M_{35}.
\end{equation}
The determinant of the $(10\times 10)$-matrix  $U_{10}$ of coefficients in
$f_{xxx},\dots, f_{ppp}$ in the functions
(\ref{nM}) is
$$
\det U_{10}=-\dfrac{3^{20}I_0^{30}}{I_1^{25}}.
$$
It does not vanish at regular points.
Since the third order invariants (\ref{nM})
are linearly independent on the fibers of the bundle
$J^3(\widetilde{\pi})\to J^2(\widetilde{\pi})$,
one can choose them as coordinates on these fibers.
The intersection  of the regular orbit and the fiber
 $\mathcal{O}_{3}\cap {\mathcal F}_3$ has zero dimension,
so the functions (\ref{nM}) generate the space of invariants of third order
and separate  regular orbits.

Invariants $M_k$ and $M_{ik}$, $i=1,2,3$,
$k=1,\dots,5$,
satisfy five syzygy relations
\begin{equation}
\label{siz1}
\begin{array}{l}
6M_{34}-6M_{25}-M_4M_5+12M_4+6M_2M_3+12M_2=0,\\[5pt]
12M_{33}-12M_{15}+24M_{21}+12M_1M_4-2M_3M_5-24M_3=0,\\[5pt]
3M_{23}-3M_{14}+M_5=0,\\[5pt]
6M_{31}-6M_{13}+2M_1M_5-3M_3^2-6M_1=0,\\[5pt]
6M_{32}-6M_{24}+2M_2M_5-3M_4^2+6M_2=0.
\end{array}
\end{equation}
The first three of them
follow from the commutation relations
(\ref{fn123}) and the Jacobi identity
for derivations $\nabla_1$, $\nabla_2$, $\nabla_3$.
The last two can be checked directly.

The situation in the orders 4 and higher is similar.
Differential invariants obtained by invariant derivations of
invariants of smaller order, are linear in fiber coordinates.
The fiber ${\cal F}_n$ of the bundle $J^n(\widetilde{\pi})\to J^{n-1}(\widetilde{\pi})$
has dimension
$\frac12(n+1)(n+2)$.
Let $N_1$, \dots, $N_{n(n+1)/2}$ be generators in the space of  invariants
of pure order $n-1$.
One can check that among the   symbols of invariants
$\nabla_iN_k$, $i=1,2,3$, $k=1,\dots, \frac12n(n+1)$,
there exist $\frac12(n+1)(n+2)$ linearly independent ones.
It follows that the corresponding invariants  generate the algebra of
invariants of  pure order $n$ and separate regular orbits.
$\Box$

\smallskip
Consider the space $\mathbb{R}^3$ with coordinates $(x,y,p)$ and
the space $\mathbb{R}^{10}$ with coordinates $(m_1,m_2,m_3,
m_4, m_5, m_{11}, m_{12}, m_{13}, m_{21}, m_{22})$.
For every function  $f$ we define the map
$$
\sigma_f: \mathbb{R}^3 \to\mathbb{R}^{10}
$$
by
$$
m_i= M_i^f, \quad m_{ik}=(\nabla_iM_k)^f,
$$
where the upper index $f$ means that the invariants are
evaluated at $f$.

Let $\Phi\in SL(3)$.
From the definition of the invariant it follows that
$$
\sigma_f \circ \Phi = \sigma_{\Phi^*f}.
$$
Hence, the image
$$
\Sigma_f = {\rm im}(\sigma_f)\subset \mathbb{R}^{10}
$$
depends only on equivalence class of $f$.

\smallskip
{\bf Definition.}
We say that the germ of  $f$ is {\it regular}
at a point $a\in \mathbb{R}^3$, if

i) 3-jets of $f$ belong to regular orbits;

ii) the germ $\sigma_f(D)$ is a germ of  smooth 3-dimensional
manifold in $\mathbb{R}^{10}$ for a domain $D\subset\mathbb{R}^3$, containing $a$;

iii) the germs of three functions of five $m_1$, $m_2$, $m_3$, $m_4$,
$m_5$  are coordinates on
$\Sigma_f$.

\smallskip
{\bf Theorem 4.}
{\it Two regular germs of functions $f$ and $\overline f$
are locally  $SL(3)$-equivalent if and only if}
\begin{equation}
\label{Sf}
\Sigma_{f}=\Sigma_{\overline{f}}.
\end{equation}

\smallskip
{\bf Proof.}
The necessity is obvious.

Assume that~(\ref{Sf}) holds.
Let us show that $f$ and $\overline{f}$ are equivalent.

Without loss of generality, we may suppose that
the germs of functions $m_1$, $m_2$, $m_3$ are local coordinates on
$\Sigma_f$.
Let
\begin{equation}
\label{KJa}
M_4^f=m_4^f(M_1, M_2, M_3),\quad M_5^f=m_5^f(M_1, M_2, M_3),\quad
M_{ik}^f=m_{ik}^f(M_1, M_2, M_3)
\end{equation}
on the submanifold $\Sigma_f$ and
$$
M_4^{\overline{f}}=m_4^{\overline{f}}(M_1, M_2, M_3),\quad
M_5^{\overline{f}}=m_5^{\overline{f}}(M_1, M_2, M_3),\quad
M_{ik}^{\overline{f}}=m_{ik}^{\overline{f}}(M_1, M_2, M_3)
$$
on $\Sigma_{{\overline{f}}}$.
The condition~(\ref{Sf}) means that
\begin{equation}
\label{KJ}
M_4^f=M_4^{\overline{f}},\quad M_5^f=M_5^{\overline{f}},\quad  M_{ik}^f=M_{ik}^{\overline{f}}.
\end{equation}

The formulas (\ref{KJa}) determine the system of 10
PDEs of third order on the function $f$.
This system is of finite type (see~\cite{Pom}).
It represents invariants $M_4$, $M_5$ and invariant derivations
$\nabla_i$ in terms of coordinates $(M_1, M_2, M_3)$.
The functions $f$ and $\overline{f}$ are solutions of this
system.
It follows from (\ref{KJ})
that invariants of all orders for $f$ and $\overline{f}$
are equal.
Since this system is of finite type, any solution is determined
by  its projection onto $J^2(\widetilde{\pi})$.
The action of $SL(3)$ on the fiber of the bundle $J^2(\widetilde{\pi}) \to
\mathbb{R}^3$ is transitive.
Consequently, $f$ and $\overline{f}$ lie in the same
orbit, hence are $SL(3)$-equivalent.
$\Box$

\smallskip
{\bf Remark.}
Note that from the formulas  (\ref{fGgp})
it follows that
the results of Section 4
can be transferred to the classification of ODEs
$y''=F(x,y,y')$ under the action of $SL(3)$
as a subgroup of pseudogroup of point transformations.
Thus, Theorems 3 and 4  also solve this classification problem.

\subsection{The case of low-dimensional group}

Now we pass to the case of two- and three-dimensional Lie groups $\cal G$ admitted by the
second order ODE.
The classification of real three-dimensional Lie algebras was performed by
Lie (see~\cite{Lie}).
We use the notations of the paper~\cite{Mah}.

Table 1. The list of non-isomorphic real three-dimensional Lie algebras.

\begin{tabular}{ll}
\hline
Algebra {\vrule width0pt height 15pt} & Non-zero commutation relations \\
\hline
$L_{3,1}$ & \\
$L_{3,2}$ &  $[X_2, X_3]=X_1$\\
$L_{3,3}$ & $[X_1, X_3]=X_1$, $[X_2, X_3]=X_1+X_2$\\
$L_{3,4}$ & $[X_1, X_3]=X_1$ \\
$L_{3,5}$ & $[X_1, X_3]=X_1$, $[X_2, X_3]=X_2$  \\
$L_{3,6}$ & $[X_1, X_3]=X_1$, $[X_2, X_3]=aX_2$, $a\ne0,1$ \\
$L_{3,7}$ & $[X_1, X_3]=aX_1-X_2$, $[X_2, X_3]=X_1+aX_2$ \\
$L_{3,8}$ & $[X_1, X_2]=X_1$, $[X_2, X_3]=X_3$, $[X_3, X_1]=-2X_2$ \\
$L_{3,9}$ & $[X_1, X_2]=X_3$, $[X_2, X_3]=X_1$, $[X_3, X_1]=X_2$ \\
\hline
\end{tabular}

\bigskip

We take the realizations of these algebras in the real plane from the
paper~\cite{Mah} (see also \cite{M-L, GKO}).
Not all algebras are presented in the following table. For the details and reasons we also refer
to~\cite{Mah}.

\bigskip
Table 2. Realizations of  three-dimensional
algebras in the real plane which are admitted by a second order
ODE ($\partial_x=\dfrac{\partial}{\partial x}$,
$\partial_y=\dfrac{\partial}{\partial y}$, $A={\rm const}$).

\begin{tabular}{lll}
\hline
Algebra {\vrule width0pt height 15pt}
& Canonical forms of generators
& Representative equations \\
\hline
$L_{3,3}^{\rm I}${\vrule width0pt height 15pt}
&
$X_1=\partial_x$, ~
$X_2=\partial_y$,
&
$y''= Ae^{-y}$
\\[3pt]
& $X_3=x\partial_x+ (x+y)\partial_y$
\\[6pt]
$L_{3,6}^{\rm I}$
&
$X_1=\partial_x$, ~
$X_2=\partial_y$, ~
&
$y''= Ay'{}^{(a-2)/(a-1)}$
\\[3pt]
& $X_3=x\partial_x+ay\partial_y$, ~
$a\ne 0,1,2$~~~~~
\\[6pt]
$L_{3,7}^{\rm I}$
&
$X_1=\partial_x$, ~
$X_2=\partial_y$,
&
$y''= A(1+y'{}^2)^{3/2}e^{a \arctan y'}$
\\[3pt]
&
$X_3=(ax+y)\partial_x+
(ay-x)\partial_y$
\\[6pt]
$L_{3,8}^{\rm I}$
&
$X_1=\partial_y$, ~
$X_2=x\partial_x+y\partial_y$,
&
$xy''= Ay'{}^3- \dfrac12 y'$
\\[3pt]
&
$X_3=2xy\partial_x+
y^2\partial_y$
\\[6pt]
$L_{3,8}^{\rm II}$
&
$X_1=\partial_y$, ~
$X_2=x\partial_x+y\partial_y$,
&
$xy''=y'+y'^3 +A(1+y'{}^2)^{3/2}$
\\[3pt]
&
$X_3=2xy\partial_x+
(y^2-x^2)\partial_y$
\\[6pt]
$L_{3,8}^{\rm III}$
&
$X_1=\partial_y$, ~
$X_2=x\partial_x+y\partial_y$,
&
$xy''=y'-y'^3 +A(1-y'{}^2)^{3/2}$
\\[3pt]
&
$X_3=2xy\partial_x+
(y^2+x^2)\partial_y$
\\[6pt]
$L_{3,9}$
&
$X_1=(1+x^2)\partial_x+xy\partial_y$,
&
$y''
= A\Bigl(
\dfrac{1+y'{}^2+(y-xy')^2}{1+x^2 +y^2}\Bigr)^{3/2}$
\\[3pt]
&
$X_2=xy\partial_x+(1+y^2)\partial_y$,
\\[3pt]
&
$X_3=y\partial_x-x\partial_y$
\\[3pt]
\hline
\end{tabular}

\bigskip

As an example we will consider in details the case of the algebra $L_{3,7}^{\rm I}$.
All other cases are similar.
Let ${\cal G}_{3,7}$ be a Lie group which Lie algebra is
$L_{3,7}^{\rm I}$.
Let us assume that
the coordinate system  $(x,y,p)$ in  $\mathbb{R}^3$
is chosen in such a way that the basis of $L_{3,7}^{\rm I}$ is
$$
X_1=\dfrac{\partial}{\partial x}, \quad
X_2=\dfrac{\partial}{\partial y}, \quad
X_3=(ax+y)\dfrac{\partial}{\partial x}+
(ay-x)\dfrac{\partial}{\partial y}.
$$

\smallskip
{\bf Proposition 3.}
{\it The action of ${\cal G}_{3,7}$ on $\widetilde{\pi}$ has one
basic invariant of order zero:}
$$
I=\dfrac{fe^{-a\arctan p}}{(1+p^2)^{3/2}}.
$$
\smallskip

Using again the Proposition 1, we find three basic $L_{3,7}^{\rm I}$-invariant
derivations:
$$
\nabla_1=\dfrac{e^{-a\arctan p}}{\sqrt{1+p^2}}\Bigl(
p\dfrac{d}{dx}-\dfrac{d}{dy}\Bigr),\quad
\nabla_2=\dfrac{e^{-a\arctan p}}{\sqrt{1+p^2}}\Bigl(
\dfrac{d}{dx}+p\dfrac{d}{dy}\Bigr),\quad
\nabla_3=(1+p^2)\dfrac{d}{dp}.
$$
They satisfy the following commutation relations
$$
[\nabla_1,\nabla_2]=0,\quad
[\nabla_1,\nabla_3]=a\nabla_1-\nabla_2,\quad
[\nabla_2,\nabla_3]=\nabla_1+a\nabla_2.
$$

Denote $I_k=\nabla_kI$ for $k=1,2,3$.
Invariants $I_{jk}=\nabla_j\nabla_kI$ are linear in second
derivatives.
The six invariants
$$
I_{11}, I_{12}, I_{13}, I_{22}, I_{23}, I_{33}
$$
have linear independent symbols and
generate the space of invariants of second order
and separate  orbits.
The determinant of the $(6\times6)$-matrix consisting of these
symbols, equals $(1+p^2)^8$ and vanishes nowhere on the fibers.

The invariants of order $\le2$ satisfy the following syzygy relations:
$$
I_{21}=I_{12},\quad I_{31}=I_{13}+I_2-aI_1,\quad I_{32}=I_{23}-I_1-aI_2.
$$

\smallskip
{\bf Theorem 5.}
{\it
Algebra of differential invariants of the
action of  ${\cal G}_{3,7}$ on $\widetilde{\pi}$
is generated by the invariant $I$ and three invariant derivations $\nabla_1$,
$\nabla_2$, $\nabla_3$.
This algebra separates  orbits.}

The proof is similar to that of Theorem 3.

\smallskip
Consider the space $\mathbb{R}^3$ with coordinates $(x,y,p)$ and
the space $\mathbb{R}^{4}$ with coordinates $(i, i_1, i_2, i_3)$.
For every function  $f$ we define the map
$$
\sigma_f: \mathbb{R}^3 \to\mathbb{R}^{4}
$$
by
$$
i=I^f, \quad i_k= I_k^f, 
$$
where the upper index $f$ means that the invariants are
evaluated at $f$.
Again, the image
$$
\Sigma_f = {\rm im}(\sigma_f)\subset \mathbb{R}^{4}
$$
depends only on equivalence class of $f$.

\smallskip
{\bf Definition.}
We say that the germ of  $f$ is {\it regular}
at a point $a\in \mathbb{R}^3$, if

i)  the germ $\sigma_f(D)$ is a germ of  smooth 3-dimensional
manifold in $\mathbb{R}^{4}$ for a domain $D\subset
\mathbb{R}^3$, containing $a$;

ii) the germs of three functions of four $i$, $i_1$, $i_2$, $i_3$ are coordinates on
$\Sigma_f$.

\smallskip
{\bf Theorem 6.}
{\it Two regular germs of functions $f$ and $\overline f$
are locally  ${\cal G}_{3,7}$-equivalent if and only if}
$$
\Sigma_{f}=\Sigma_{\overline{f}}.
$$

The proof is similar to that of Theorem 4.


\medskip

In what follows we present the basic invariant $I$, three invariant derivations and their
commutation relations for the rest of the algebras from the Table 2.

1. Algebra $L_{3,3}^{\rm I}$.

Basic invariant
$$
I=fe^p.
$$

Invariant derivations:
$$
\nabla_1=e^p\dfrac{d}{dx}+pe^p\dfrac{d}{dy},\quad
\nabla_2=e^p\dfrac{d}{dy},\quad
\nabla_3=\dfrac{d}{dp}.
$$

Commutation relations:
$$
[\nabla_1,\nabla_2]=0,\quad
[\nabla_1,\nabla_3]=-\nabla_1-\nabla_2,\quad
[\nabla_2,\nabla_3]= -\nabla_2.
$$



2. Algebra $L_{3,6}^{\rm I}$.

Basic invariant
$$
I=fp^{(2-a)/(a-1)}.
$$

Invariant derivations:
$$
\nabla_1=p^{1/(a-1)}\dfrac{d}{dx},\quad
\nabla_2=p^{a/(a-1)}\dfrac{d}{dy},\quad
\nabla_3=p\dfrac{d}{dp}.
$$

Commutation relations:
$$
[\nabla_1,\nabla_2]=0,\quad
[\nabla_1,\nabla_3]=\frac{1}{1-a}\nabla_1,\quad
[\nabla_2,\nabla_3]= \frac{a}{1-a}\nabla_2.
$$



3. Algebra
$L_{3,8}^{\rm I}$.

Basic invariant
$$
I=\dfrac{p+2xf}{p^3}.
$$

Invariant derivations:
$$
\nabla_1=x\dfrac{d}{dx}-p\dfrac{d}{dp},\quad
\nabla_2=\dfrac xp\dfrac{d}{dx}+x\dfrac{d}{dy}-\dfrac12\dfrac{d}{dp},\quad
\nabla_3=p^2\dfrac{d}{dp}.
$$

Commutation relations:
$$
[\nabla_1,\nabla_2]=\nabla_2,\quad
[\nabla_1,\nabla_3]=-\nabla_3,\quad
[\nabla_2,\nabla_3]=\nabla_1.
$$

4. Algebra
$L_{3,8}^{\rm II}$.

Basic invariant
$$
I=\dfrac{-fx+p+p^3}{(1+p^2)^{3/2}}.
$$

Invariant derivations:
\begin{multline*}
\nabla_1=\dfrac{x}{\sqrt{1+p^2}}\dfrac{d}{dx}+\dfrac{xp}{\sqrt{1+p^2}}\dfrac{d}{dy}
+p\sqrt{1+p^2}\dfrac{d}{dp},\\[7pt]
\nabla_2=-\dfrac{xp}{\sqrt{1+p^2}}\dfrac{d}{dx}+\dfrac{x}{\sqrt{1+p^2}}\dfrac{d}{dy}
+\sqrt{1+p^2}\dfrac{d}{dp},\quad
\nabla_3=(1+p^2)\dfrac{d}{dp}.
\end{multline*}

Commutation relations:
$$
[\nabla_1,\nabla_2]=-\nabla_3,\quad
[\nabla_1,\nabla_3]=-\nabla_2,\quad
[\nabla_2,\nabla_3]=\nabla_1.
$$

5. Algebra
$L_{3,8}^{\rm III}$.

For this algebra the domain is divided into two parts --- where
$|p|<1$ and $|p|>1$.

For $|p|<1$ the basic
invariant is
$$
I=\dfrac{fx-p+p^3}{(1-p^2)^{3/2}}
$$
and the invariant derivations are
\begin{multline*}
\nabla_1=\dfrac{x}{\sqrt{1-p^2}}\dfrac{d}{dx}+\dfrac{xp}{\sqrt{1-p^2}}\dfrac{d}{dy}
+p\sqrt{1-p^2}\dfrac{d}{dp},\\[7pt]
\nabla_2=\dfrac{xp}{\sqrt{1-p^2}}\dfrac{d}{dx}+\dfrac{x}{\sqrt{1-p^2}}\dfrac{d}{dy}
+\sqrt{1-p^2}\dfrac{d}{dp},\quad
\nabla_3=(1-p^2)\dfrac{d}{dp}.
\end{multline*}
Their commutation relations are:
$$
[\nabla_1,\nabla_2]=-\nabla_3,\quad
[\nabla_1,\nabla_3]=-\nabla_2,\quad
[\nabla_2,\nabla_3]=-\nabla_1.
$$

For $|p|>1$
the basic invariant is
$$
I=\dfrac{fx-p+p^3}{(p^2-1)^{3/2}}
$$
and the invariant derivations are
\begin{multline*}
\nabla_1=\dfrac{x}{\sqrt{p^2-1}}\dfrac{d}{dx}+\dfrac{xp}{\sqrt{p^2-1}}\dfrac{d}{dy}
-p\sqrt{p^2-1}\dfrac{d}{dp},\\[7pt]
\nabla_2=\dfrac{xp}{\sqrt{p^2-1}}\dfrac{d}{dx}+\dfrac{x}{\sqrt{p^2-1}}\dfrac{d}{dy}
-\sqrt{p^2-1}\dfrac{d}{dp},\quad
\nabla_3=(p^2-1)\dfrac{d}{dp}.
\end{multline*}
Their commutation relations are:
$$
[\nabla_1,\nabla_2]=-\nabla_3,\quad
[\nabla_1,\nabla_3]=\nabla_2,\quad
[\nabla_2,\nabla_3]=-\nabla_1.
$$

6. Algebra $L_{3,9}$.


This case is slightly more complicated.
We have three relative invariants of order zero:
$$
I_0=x^2+y^2+1,\quad I_1=p^2(1+x^2)-2xyp+y^2+1=p^2+(px-y)^2+1, \quad I_2=f.
$$
Using them, we construct the basic invariant
$$
I=\dfrac{f\cdot I_0^{3/2}}{I_1^{3/2}}.
$$
Note that $I_0>0$, $I_1>0$ everywhere.

The invariant derivations are:
\begin{multline*}
\nabla_1=\dfrac{\sqrt{I_0}}{\sqrt{I_1}}\Bigl(
(p(x^2+1)-xy)\dfrac{d}{dx}+(xyp-(1+y^2))\dfrac{d}{dy}\Bigr)
-\dfrac{(x+yp)\sqrt{I_1}}{\sqrt{I_0}}\dfrac{d}{dp},\\[7pt]
\nabla_2=\dfrac{I_0}{\sqrt{I_1}}\Bigl(
\dfrac{d}{dx}+p\dfrac{d}{dy}\Bigr),\quad
\nabla_3=\dfrac{I_1}{\sqrt{I_0}}\dfrac{d}{dp}.
\end{multline*}
Their commutation relations are:
$$
[\nabla_1,\nabla_2]=\nabla_3,\quad
[\nabla_1,\nabla_3]=-\nabla_2,\quad
[\nabla_2,\nabla_3]=\nabla_1.
$$

In all of the cases above the equivalence theorem is formulated in the same manner
as the Theorem 6.

\bigskip
Now we consider the case when associated second order equation admits the
two-dimensional Lie algebra.
It is well-known (see, e.g.,~\cite{Ibr09, Mah}) that there are four following
realizations of two-dimensional Lie algebra in $\mathbb{R}^2$.

\bigskip
Table 3. Realizations of  two-dimensional
algebras in the real plane which are admitted by a second order
ODE.

\medskip
\begin{tabular}{lll}
\hline
Algebra {\vrule width0pt height 15pt}
& Canonical forms of generators\mbox{\qquad}
& Representative equations \\
\hline
$L_{2,1}^{\rm I}${\vrule width0pt height 15pt}
&
$X_1=\partial_x$, ~
$X_2=\partial_y$
&
$y''= f(y')$
\\[3pt]
$L_{2,1}^{\rm II}$
&
$X_1=\partial_x$, ~
$X_2=x\partial_y$
&
$y''= f(x)$
\\[3pt]
$L_{2,2}^{\rm I}$
&
$X_1=\partial_y$, ~
$X_2=x\partial_x+y\partial_y$
&
$xy''=f(y')$
\\[3pt]
$L_{2,2}^{\rm II}$
&
$X_1=\partial_y$, ~
$X_2=y\partial_y$
&
$y''= y'f(x)$
\\[3pt]
\hline
\end{tabular}

\bigskip

Lie showed that the algebras $L_{2,1}^{\rm II}$ and $L_{2,2}^{\rm II}$
result in linearization
of the underlying second-order equation.
So, it remains to consider the other two realizations.

For the algebra $L_{2,1}^{\rm I}$ we have two basic invariants of order zero
$$
I=p, \quad J=f
$$
and three invariant derivations
$$
\nabla_1=\dfrac{d}{dx},\quad
\nabla_2=\dfrac{d}{dy},\quad
\nabla_3=\dfrac{d}{dp}.
$$
They all obviously commute.

For the algebra $L_{2,2}^{\rm I}$
the basic invariants are
$$
I=p, \quad J=xf
$$
and the invariant derivations are
$$
\nabla_1=x\dfrac{d}{dx},\quad
\nabla_2=x\dfrac{d}{dy},\quad
\nabla_3=\dfrac{d}{dp}.
$$
They satisfy the following commutation relations:
$$
[\nabla_1,\nabla_2]=\nabla_2,\quad
[\nabla_1,\nabla_3]=0,\quad
[\nabla_2,\nabla_3]=0.
$$

In both cases the map
$$
\sigma_f: \mathbb{R}^3=\{(x,y,p)\} \to
\mathbb{R}^{5}=\{(i, j,  j_1, j_2, j_3)\}
$$
is defined by
$$
i=I^f, \quad j=J^f,\quad j_k= (\nabla_kJ)^f.
$$
Again, we define $\Sigma_f = {\rm im}(\sigma_f)\subset \mathbb{R}^{5}$.

The equivalence theorem is formulated in the same manner as Theorem 6:
two regular germs of functions $f$ and $\overline f$
are locally  equivalent if and only if
$\Sigma_{f}=\Sigma_{\overline{f}}$.

\smallskip
{\bf Acknowledgement.}
The author would like to express his deep gratitude to Professor
Valentin Lychagin for his guidance
throughout the work and valuable comments.

\vskip1cm

\begin{flushleft}
\rm
Institute of Mathematics and Mechanics,\\
Kazan (Volga Region) Federal University,\\
 Kazan, Russia\\[10pt]
{\it E-mail:}~ {\tt vshjr@yandex.ru}, {\tt 1Vadim.Shurygin@kpfu.ru}.
\end{flushleft}

\end{document}